\documentclass[12pt]{article}
\usepackage{amssymb,amsmath}
\usepackage{graphicx}
\usepackage{amsbsy}
\usepackage{enumerate}

\newtheorem{theorem}{Theorem}[section]
\newtheorem{example}[theorem]{Example}
\newtheorem{proposition}[theorem]{Proposition}

\newtheorem{definition}[theorem]{Definition}
\newtheorem{remark}[theorem]{Remark}
\newtheorem{step}{Step}
\numberwithin{equation}{section}

\begin{document}

\begin{center}{\Large\sc
Caffarelli-Kohn-Nirenberg inequality\\
on metric measure spaces with applications
}\\
\vspace{0.5cm}
 {\large  Alexandru Krist\'aly\\
 {\normalsize Department of Economics, Babe\c s-Bolyai University, 400591 Cluj-Napoca,
Romania\\
Email addresses: alexandrukristaly@yahoo.com;  alexandru.kristaly@econ.ubbcluj.ro}}\\
\vspace{0.5cm} {\large  Shin-ichi Ohta}\\
 {\normalsize Department of Mathematics, Kyoto University, Kyoto 606-8502, Japan\\
 Email address: sohta@math.kyoto-u.ac.jp}
\end{center}

\vspace{0.1cm}

\begin{center}
{\small\it  Dedicated to Professor Vicen\c tiu R\u adulescu on the
occasion of his $55^{th}$
 birthday}
 \end{center}

 \vspace{0.1cm}

\begin{abstract}
\noindent We prove that if a metric measure space satisfies the
volume doubling condition and the Caffarelli-Kohn-Nirenberg
inequality with the same exponent $n \ge 3$, then it has exactly the
$n$-dimensional volume growth. As an application, if an
$n$-dimensional Finsler manifold of non-negative $n$-Ricci curvature
satisfies the Caffarelli-Kohn-Nirenberg inequality with the sharp
constant, then its flag curvature is identically zero. In the
particular case of Berwald spaces, such a space is necessarily
isometric to a Minkowski space.
\bigskip

\noindent {\it Keywords}: Caffarelli-Kohn-Nirenberg inequality;
metric measure spaces; volume doubling condition; Finsler manifolds.

\noindent {\it MSC}: 35R06, 53C60, 58J60.
\end{abstract}

\section{Introduction and statement of main results}

Let $a\in [0,1)$ be a parameter, $n\geq 3$ be an integer, and put
$p={2n}/({n-2+2a})$. In the theory of Sobolev inequalities, a
central role is played by the famous \emph{Caffarelli-Kohn-Nirenberg
inequality} (see \cite{CKN-cikk}) which states that
\[
\left(\int_{\mathbb R^n} \frac{|u|^p}{|x|^{ap}}
\,dx\right)^\frac{1}{p} \leq K_a \left(\int_{\mathbb R^n} |D u|^2
\,dx\right)^\frac{1}{2}\quad \text{for all}\ u\in C_0^\infty(\mathbb
R^n),
\]
where
\begin{equation}\label{ckn-eukl}
K_a:=\left(\frac{1}{(n-2)(n-ap)}\right)^\frac{1}{2}
\left(\frac{(2-ap)\Gamma((2n-2ap)/(2-ap))}{n\omega_n\Gamma^2((n-ap)/(2-ap))}\right)^\frac{2-ap}{2n-2ap}
\end{equation}
is the optimal constant (see Lieb \cite{Lieb}),
$\omega_n:=\pi^{n/2}/\Gamma(n/2 +1)$ being the volume of the unit
ball in $\mathbb R^n$. Moreover, a family of extremals is given by
\begin{equation}\label{minim-eukl}
u_\lambda(x)={\left(\lambda+|x|^{2-ap}\right)^\frac{2-n}{2-ap}},
\quad \lambda>0.
\end{equation}
The optimal constant and extremals for $a=0$ have been established
by Aubin~\cite{Aubin} and Talenti~\cite{Talenti} in which case the
above inequality reduces to the standard Sobolev inequality; see
Chou and Chu~\cite{CC-London} for the most general case.
Furthermore, various versions of the Caffarelli-Kohn-Nirenberg
inequality have been treated also on Riemannian manifolds and
Orlicz-Sobolev spaces (see, e.g., do Carmo and
Xia~\cite{doCarmo-Xia}).

The main objective of the present paper is to investigate the
Caffarelli-Kohn-Nirenberg inequality in the context of metric
measure spaces. As applications, we provide novel rigidity results
for Finsler manifolds by means of the \emph{sharp}
Caffarelli-Kohn-Nirenberg inequality.

In order to state the main result of the paper, we fix the numbers
$a,n$ and $p$ as above. Let $(X,d)$ be a metric space and $\mu$ be a
Borel measure on $X$ such that $0<\mu(U)<\infty$ for any nonempty
bounded open set $U \subset X$. For some element $x_0\in X$ and
constant $C>0$, we consider the Caffarelli-Kohn-Nirenberg inequality
on $(X,d,\mu)$ of the form
\[
\left(\int_X \frac{|u(x)|^p}{d(x_0,x)^{ap}}
\,d\mu(x)\right)^\frac{1}{p} \leq C \left(\int_X |D u|(x)^2
\,d\mu(x)\right)^\frac{1}{2}\ \text{for all}\ u\in {\rm Lip}_0(X).
\eqno{({\bf CKN})_C^{x_0}}
\]
Hereafter, ${\rm Lip}_0(X)$ is the space of Lipschitz functions with
compact support on $X$, while
\[ |Du|(x):=\limsup_{y\to x} \frac{|u(y)-u(x)|}{d(x,y)} \]
is the local Lipschitz constant of $u$ at $x\in X$. The function
$x\longmapsto |Du|(x)$ is Borel measurable for $u\in {\rm
Lip}_0(X)$. For instance, any bi-Lipschitz deformation of the
Euclidean space $\mathbb R^n$ satisfies $({\bf CKN})_C^{x_0}$ with
some $C \ge K_a$.

For some fixed elements $C_0\geq 1$ and $x_0\in X$, we introduce the
following hypotheses on the behavior of $\mu$:
\begin{description}
\item[{\rm $({\bf VD})^n_{C_0}$}]
 $\displaystyle\frac{\mu(B(x,R))}{\mu(B(x,r))}\leq {C_0}
 \left(\frac{R}{r}\right)^n$ for all  $x\in X$ and $0<r<R$;
\item[{\rm $({\bf AR})^n_{x_0}$}]
 $\displaystyle\liminf_{r\to 0}\frac{\mu(B(x_0,r))}{\mu_E(\mathbb B_n(r))}=1$.
\end{description}
As usual, $B(x,r):=\{y\in X:d(x,y)<r\}$, $\mathbb B_n(r):=\{x\in
\mathbb R^n:|x|<r\}$, and $\mu_E$ is the $n$-dimensional Lebesgue
measure.

The main result of the paper can be stated as follows.

\begin{theorem}\label{theorem-1-fo}
Let $a\in [0,1)$, $n \geq 3$, $p={2n}/({n-2+2a})$, $x_0\in X$,
$C\geq K_a$, and $C_0\geq 1$. Assume that the
Caffarelli-Kohn-Nirenberg inequality
$({\bf CKN})_C^{x_0}$ holds on a proper 
metric measure space $(X,d,\mu)$, and the hypotheses $({\bf
VD})^n_{C_0}$ and $({\bf AR})^n_{x_0}$ are verified. Then, for every
$x\in X$ and $\rho>0$, we have
\begin{equation}\label{volume-estimate}
\mu(B(x,\rho)) \ge
 C_0^{-1}(C^{-1}K_a)^\frac{n}{1-a}\mu_E(\mathbb B_n(\rho)).
\end{equation}
In particular, $(X,d,\mu)$ has the $n$-dimensional volume growth
\[
C_0^{-1}(C^{-1}K_a)^\frac{n}{1-a}\omega_n \rho^n
 \le \mu(B(x_0,\rho)) \le C_0 \omega_n \rho^n
 \quad \text{for all}\ \rho>0.
\]
\end{theorem}

This theorem extends do Carmo and Xia's result
\cite[Theorem~1.1]{doCarmo-Xia} on Riemannian manifolds of
non-negative Ricci curvature in two respects.
Theorem~\ref{theorem-1-fo} is concerned with general metric measure
spaces, and assumes only the volume growth condition $({\bf
VD})^n_{C_0}$ instead of the curvature bound. Before discussing
applications, let us give several remarks on the hypotheses and the
conclusions of the theorem.

\begin{remark}\rm
(a) We remark that $({\bf CKN})_C^{x_0}$ ensures that $(X,d)$ is
unbounded (equivalently, non-compact). Indeed, if $(X,d)$ is
bounded, then $u+c$ with $c \to \infty$ violates the validity of
$({\bf CKN})_C^{x_0}$.

(b) If $(X,d,\mu)$ satisfies the \emph{volume doubling condition}:
\[ \mu(B(x,2r)) \le \Lambda \mu(B(x,r)) \quad \text{for some}\ \Lambda \ge 1\
 \text{and all}\ x \in X,\ r>0, \]
then we easily see that $({\bf VD})^n_{C_0}$ is satisfied (with,
e.g., $n \ge \log_2 \Lambda$ and $C_0=\Lambda$). Thus $({\bf
VD})^n_{C_0}$ can be interpreted as the volume doubling condition
with the explicit exponent $n$. One can also regard $({\bf
VD})^n_{C_0}$ as a generalization of the Bishop-Gromov volume growth
estimate (of non-negative Ricci curvature).

(c) Note that, on the one hand, $({\bf VD})^n_{C_0}$ implies that
the Hausdorff dimension $\dim_H X$ of $(X,d)$ is at most $n$. On the
other hand, since $\limsup_{r\to 0}\mu(B(x_0,r))/\mu_E(\mathbb
B_n(r)) \le C_0$ by $({\bf VD})^n_{C_0}$ and $({\bf AR})^n_{x_0}$,
we have the \emph{Ahlfors $n$-regularity} at $x_0$ in the sense that
$\Omega^{-1}r^n \le \mu(B(x_0,r)) \le \Omega r^n$ for some $\Omega
\ge 1$ and small $r>0$, thus we have $\dim_H X=n$. (See
\cite{Heinonen} for the importance of the volume doubling condition
and the Ahlfors regularity in analysis on metric measure spaces.) We
also remark that the constant $1$ was chosen as the RHS of $({\bf
AR})^n_{x_0}$ merely for simplicity. Since
$\Omega_{x_0}:=\liminf_{r\to 0}\mu(B(x_0,r))/r^n$ is positive by
$({\bf VD})^n_{C_0}$, one can normalize $\mu$ so as to satisfy
$({\bf AR})^n_{x_0}$ whenever $\Omega_{x_0}$ is bounded.

(d) The volume growth estimate \eqref{volume-estimate} shows that,
for instance, the cylinder $\mathbb S^{n-1} \times \mathbb R$ can
not support $({\bf CKN})_C^x$ for \emph{any} $x$ and $C$.

(e) The use of the measure $\mu_E$ as a comparing one to $\mu$ in
the hypotheses $({\bf VD})^n_{C_0}$ and $({\bf AR})^n_{x_0}$ comes
from the fact that the number $K_a$ is optimal and the functions
from (\ref{minim-eukl}) are minimizers in the {\it Euclidean}
Caffarelli-Kohn-Nirenberg inequality. Therefore, if a sharp
Caffarelli-Kohn-Nirenberg inequality holds in a generic metric
measure space $(X_0,d_0,\mu_0)$, knowing the optimal constant
$K_0>0$ and assuming that the class of extremals is formally the
same as (\ref{minim-eukl}) with $d_0(x_0,x)$ instead of $|x|$, then
one can prove a similar statement to Theorem~\ref{theorem-1-fo} by
replacing $\mu_E$ and $K_a$ with $\mu_0$ and $K_0$, respectively.
\end{remark}

We point out that, on (absolutely homogeneous for simplicity)
Finsler manifolds with non-negative \emph{$n$-Ricci curvature},
$({\bf VD})^n_{C_0}$ holds with $C_0=1$ (see
Shen~\cite{Shen-volume}, Ohta~\cite{Ohta-Fint} and
Theorem~\ref{volume-comp-fo} below). In particular, from
Theorem~\ref{theorem-1-fo}, important rigidity results can be
deduced in the context of Finsler manifolds when the sharp
Caffarelli-Kohn-Nirenberg inequality holds (for precise notions, see
Section~\ref{sect-3}). We state two such results.

\begin{theorem}\label{theorem-Finsler}
Let $a\in [0,1)$, $n \ge 3$, $p={2n}/({n-2+2a})$, and $(M,F)$ be a
complete $n$-dimensional Finsler manifold. Fix a positive smooth
measure $\mu$ on $M$ and assume that the $n$-Ricci curvature
$\mathrm{Ric}_n$ of $(M,F,\mu)$ is non-negative, the sharp
Caffarelli-Kohn-Nirenberg inequality $({\bf CKN})_{K_a}^{x_0}$ holds
for some $x_0\in M$, and in addition $\lim_{r \to
0}\mu(B(x_0,r))/\omega_n r^n=1$. Then the flag curvature of $(M,F)$
is identically zero.
\end{theorem}

\begin{theorem}\label{theorem-Berwald}
Let $a\in [0,1)$, $n \ge 3$, $p={2n}/({n-2+2a})$, and $(M,F)$ be a
complete $n$-dimensional Berwald space with non-negative Ricci
curvature. If for some $x_0\in M$ and the $n$-dimensional Hausdorff
measure of $(M,F)$ the sharp Caffarelli-Kohn-Nirenberg inequality
$({\bf CKN})_{K_a}^{x_0}$ holds, then $(M,F)$ is isometric to a
Minkowski space.
\end{theorem}

\begin{remark}\rm
(a) By using anisotropic symmetrization arguments, we prove in
Section~\ref{sect-3} that the sharp Caffarelli-Kohn-Nirenberg
inequality $({\bf CKN})_{K_a}^{x_0}$ holds on every Minkowski space
$(\mathbb R^n,F)$ (Proposition~\ref{prop-1-Minko}). In this manner,
Theorem~\ref{theorem-Berwald} delimits Minkowski spaces as the
optimal geometric framework where $({\bf CKN})_{K_a}^{x_0}$ holds
within the class of Berwald spaces with non-negative Ricci
curvature.

(b) Riemannian manifolds and (locally) Minkowski spaces are Berwald
spaces. Therefore Theorem~\ref{theorem-Berwald} extends do Carmo and
Xia's result \cite[Corollary~1.2]{doCarmo-Xia} in the Riemannian
context. In fact, some constructions in the proof of
Theorem~\ref{theorem-1-fo} are inspired from \cite{doCarmo-Xia}.

(c) The assumption $n \geq 3$ in Theorem~\ref{theorem-Berwald} is
essential not only in the definition of $p={2n}/({n-2+2a})$ but also
for the structure of the Berwald space. Indeed, Szab\'o's rigidity
result states that any Berwald {\it surface} is either a locally
Minkowski space or a Riemannian surface, see Szab\'o~\cite{Szabo}
and Bao, Chern and Shen~\cite[Theorem~10.6.2]{BCS}.
\end{remark}

The paper is constructed as follows. In Section~\ref{sect-2}, we
prove Theorem~\ref{theorem-1-fo}. In Section~\ref{sect-3}, we first
recall some basic notions and results from Finsler geometry, and
then complete the proof of Theorems~\ref{theorem-Finsler} and
\ref{theorem-Berwald}.

\section{Proof of Theorem~\ref{theorem-1-fo}}\label{sect-2}

We divide the proof into five steps.

\begin{step}\rm
We first derive an important ODE from the extremals
\eqref{minim-eukl} in the Euclidean case. Since $u_\lambda(x)=
{\left(\lambda+|x|^{2-ap}\right)^{\frac{2-n}{2-ap}}}$ is a minimizer
in the Euclidean Caffarelli-Kohn-Nirenberg inequality $({\bf
CKN})_{K_a}^{x_0}$, the following integral identity holds  for every
$\lambda>0$:
\begin{equation}\label{identity-1}
    \left(\int_{\mathbb R^n}
    \frac{\left(\lambda+|x|^{2-ap}\right)^{\frac{(2-n)p}{2-ap}}}{|x|^{ap}} \,d\mu_E(x)\right)^\frac{2}{p}
    =K_a^{2}(n-2)^2
    \int_{\mathbb R^n}\frac{\left(\lambda+|x|^{2-ap}\right)^{\frac{2(ap-n)}{2-ap}}}{|x|^{2ap-2}} \,d\mu_E(x).
\end{equation}
Observe that
$\frac{(2-n)p}{2-ap}=\frac{2(ap-n)}{2-ap}=\frac{n}{a-1}<0$, in
particular, $ap<2$. We introduce the auxiliary function
$Q_E:(0,\infty)\longrightarrow \mathbb R$ defined by
\[
Q_E(\lambda):=\frac{1-a}{n-1+a} \int_{\mathbb
R^n}\frac{\left(\lambda+|x|^{2-ap}\right)^{\frac{n-1+a}{a-1}}}{|x|^{ap}}
 \,d\mu_E(x).
\]
Then the identity \eqref{identity-1} reduces to, provided that $Q_E$
is well-defined,
\begin{equation}\label{minko-1-azonossag}
(-Q_E'(\lambda))^\frac{2}{p}=K_a^{2}(n-2)^2\left(\frac{n-1+a}{1-a}Q_E(\lambda)+\lambda
Q_E'(\lambda)\right),\quad \lambda>0.
\end{equation}
To see that $Q_E$ is well-defined, we obtain from the layer cake
representation of functions and a change of variables as $t=(\lambda
+\rho^{2-ap})^{\frac{n-1+a}{a-1}} \rho^{-ap}$ that
\begin{align*}
Q_E(\lambda) &= \frac{1-a}{n-1+a} \int_0^\infty
 \mu_E\left\{x\in \mathbb R^n: \frac{\left(\lambda+|x|^{2-ap}\right)^{\frac{n-1+a}{a-1}}}{|x|^{ap}}>t\right\} dt \\
 &= \frac{1-a}{n-1+a}\int_0^\infty
 \mu_E\left\{x\in \mathbb R^n:|x|<\rho
   \right\} f(\lambda,\rho) \,d\rho,
\end{align*}
where $f:(0,\infty)^2\longrightarrow \mathbb R$ is given by
\begin{align}
 f(\lambda,\rho)
 &=\frac{n-1+a}{1-a} \frac{(\lambda+\rho^{2-ap})^{\frac{n}{a-1}}}{\rho^{ap}} (2-ap) \rho^{1-ap}
 +ap \frac{(\lambda+\rho^{2-ap})^{\frac{n-1+a}{a-1}}}{\rho^{ap+1}} \nonumber\\
 &=\frac{(\lambda+\rho^{2-ap})^{\frac{n}{a-1}}}{\rho^{ap+1}}
 \left\{ \rho^{2-ap}\left(\frac{n-1+a}{1-a}(2-ap)+ap\right)+ap\lambda\right\}.
 \label{f-definic}
\end{align}
Hence we have
\begin{equation}\label{Q-definicio-3}
 Q_E(\lambda)= \frac{1-a}{n-1+a}\int_0^\infty \mu_E(\mathbb B_n(\rho)) f(\lambda,\rho) \,d\rho.
\end{equation}
An elementary calculus shows that the improper integral in
\eqref{Q-definicio-3} converges, thus $Q_E$ is well-defined.
\end{step}

\begin{step}\rm
Switching to the metric measure setting as in
Theorem~\ref{theorem-1-fo}, we first observe that the hypotheses
$({\bf VD})^n_{C_0}$ and $({\bf AR})^n_{x_0}$ yield
\begin{equation}\label{two-measures}
\mu(B(x_0,\rho))\leq C_0 \mu_E(\mathbb B_n(\rho)) \quad \text{for
every}\ \rho>0.
\end{equation}
Let us consider for each $\lambda>0$ the sequence of functions
$u_{\lambda,k}:X\longrightarrow \mathbb R$, $k\in \mathbb N$,
defined by
\[ u_{\lambda,k}(x)
 :=\max\{ 0,\min\{0,k-d(x_0,x)\}+1\}
 \left(\lambda+\max\left\{ d(x_0,x),k^{-1}\right\}^{2-ap}\right)^{\frac{2-n}{2-ap}}. \]
Since $(X,d)$ is proper, ${\rm supp}(u_{\lambda,k})=\{x\in
X:d(x_0,x)\leq k+1\}$ is compact. Therefore we have
$u_{\lambda,k}\in {\rm Lip}_0(X)$ for every $\lambda>0$ and $k\in
\mathbb N$. We set
\[ \tilde u_\lambda(x)
 :=\lim_{k\to \infty}u_{\lambda,k}(x)= {\left(\lambda+d(x_0,x)^{2-ap}\right)^\frac{2-n}{2-ap}}. \]
Since the functions $u_{\lambda,k}$ verify the
Caffarelli-Kohn-Nirenberg inequality $({\bf CKN})_C^{x_0}$, a simple
approximation procedure based on (\ref{two-measures}) shows that
$\tilde u_{\lambda}$ verifies $({\bf CKN})_C^{x_0}$ as well.
Consequently, we can apply $({\bf CKN})_C^{x_0}$ to $\tilde
u_\lambda$. In particular, by exploiting a chain rule for the local
Lipschitz constant and the fact that $x\longmapsto d(x_0,x)$ is
$1$-Lipschitz (thus $|Dd(x_0,\cdot)|(x)\leq 1$ for all $x$), we
obtain
\begin{equation}\label{ssssok} \left(\int_X
\frac{{\left(\lambda+d(x_0,x)^{2-ap}\right)^\frac{(2-n)p}{2-ap}}}{d(x_0,x)^{ap}}
\,d\mu(x)\right)^\frac{2}{p}
 \leq C^2(n-2)^2 \int_X
 \frac{\left(\lambda+d(x_0,x)^{2-ap}\right)^\frac{2(ap-n)}{2-ap}}{d(x_0,x)^{2ap-2}} \,d\mu(x).
\end{equation}
We shall rewrite \eqref{ssssok} by means of  the function $\tilde
Q:(0,\infty)\longrightarrow \mathbb R$ defined by
\[ \tilde Q(\lambda) :=\frac{1-a}{n-1+a} \int_{X}
 \frac{\left(\lambda+d(x_0,x)^{2-ap}\right)^{\frac{n-1+a}{a-1}}}{d(x_0,x)^{ap}} \,d\mu(x). \]
Before to do that, we claim that $\tilde Q$ is well-defined. Again,
by the layer cake representation of functions, one has
\[ \tilde Q(\lambda)=\frac{1-a}{n-1+a}\int_0^\infty
 \mu\left\{x\in X:
 \frac{\left(\lambda+d(x_0,x)^{2-ap}\right)^{\frac{n-1+a}{a-1}}}{d(x_0,x)^{ap}}>t\right\} dt. \]
By taking into account that $\mathrm{diam}\, X=\infty$, similarly to
the previous step, we have
\begin{equation}\label{Q-tilde-def}
\tilde Q(\lambda)=\frac{1-a}{n-1+a}\int_0^\infty
\mu(B(x_0,\rho))f(\lambda,\rho) \,d\rho.
\end{equation}
In particular, from \eqref{two-measures} and (\ref{Q-definicio-3}),
for every $\lambda>0$ we obtain
\[ 0<\tilde Q(\lambda) \leq
 \frac{C_0(1-a)}{n-1+a}\int_0^\infty \mu_E(\mathbb B_n(\rho))
 f(\lambda,\rho) \,d\rho= C_0 Q_E(\lambda), \]
which concludes the claim. Now, similarly to
\eqref{minko-1-azonossag}, we can transform the relation
\eqref{ssssok} via $\tilde Q$ into the inequality
\begin{equation}\label{minko-2-egyenlotleseg}
(-\tilde Q'(\lambda))^\frac{2}{p}\leq
C^2(n-2)^2\left(\frac{n-1+a}{1-a}\tilde Q(\lambda)+\lambda \tilde
Q'(\lambda)\right),\quad \lambda>0.
\end{equation}

Inspired from \eqref{minko-1-azonossag} and
\eqref{minko-2-egyenlotleseg}, we consider the ODE
\begin{equation}\label{minko-3-egyenloseg}
(- q'(\lambda))^\frac{2}{p}= C^2(n-2)^2\left(\frac{n-1+a}{1-a}
q(\lambda)+\lambda  q'(\lambda)\right),\quad \lambda>0.
\end{equation}
On account of \eqref{minko-1-azonossag}, one can observe that
\eqref{minko-3-egyenloseg} has the particular solution of the form
\[
 q(\lambda)=(C^{-1}K_a)^\frac{n}{1-a}Q_E(\lambda).
\]
\end{step}

\begin{step}\rm
We shall show that, for every $\lambda>0$,
\begin{equation}\label{g-kozott-egyenlotlenseg}
    \tilde Q(\lambda)\geq q(\lambda).
\end{equation}
Suppose $C>K_a$ without loss of generality. The proof of
\eqref{g-kozott-egyenlotlenseg} requires a local (near zero) and a
global treatment of the quotient $\tilde Q/q$. First, due to the
hypothesis $({\bf AR})^n_{x_0}$, for every $\varepsilon>0$ there
exists $\rho_\varepsilon>0$ such that $\mu(B(x_0,\rho))\geq
(1-\varepsilon)\mu_E(\mathbb B_n(\rho))$ for all $\rho\in
[0,\rho_\varepsilon]$. Therefore, by \eqref{Q-tilde-def} and
changing the variables as $\rho=\lambda^\frac{1}{2-ap}t$, it turns
out that
\begin{align*}
\tilde Q(\lambda)
 &\ge  \frac{1-a}{n-1+a}(1-\varepsilon)
 \int_0^{\rho_\varepsilon} \mu_E(\mathbb B_n(\rho))f(\lambda,\rho) \,d\rho \\
  &= \frac{1-a}{n-1+a}(1-\varepsilon)\lambda^\frac{n-2+2a}{2(a-1)}
  \int_0^{\rho_\varepsilon \lambda^\frac{1}{ap-2}} \mu_E(\mathbb B_n(t))f(1,t) \,dt.
\end{align*}
A similar argument gives from \eqref{Q-definicio-3} that
\begin{equation}\label{Q_E}
Q_E(\lambda)= \frac{1-a}{n-1+a}\lambda^\frac{n-2+2a}{2(a-1)}
 \int_0^\infty \mu_E(\mathbb B_n(t)) f(1,t) \,dt.
\end{equation}
The above relations and the fact $ap-2<0$ lead to
\[ \liminf_{\lambda\to 0}\frac{\tilde Q(\lambda)}{q(\lambda)}
 =(CK_a^{-1})^\frac{n}{1-a}\liminf_{\lambda\to 0}\frac{\tilde Q(\lambda)}{Q_E(\lambda)}
 \geq (CK_a^{-1})^\frac{n}{1-a}(1-\varepsilon). \]
Since $\varepsilon>0$ is arbitrarily small, we obtain
\[
    \liminf_{\lambda\to 0}\frac{\tilde Q(\lambda)}{q(\lambda)}\geq
    (CK_a^{-1})^\frac{n}{1-a}> 1,
\]
concluding the study of the quotient $\tilde Q/q$ near the origin.

Now, arguing by contradiction, we assume that there exists $\tilde
\lambda>0$ such that $\tilde Q(\tilde \lambda)<q(\tilde \lambda)$.
By the continuity of the functions $\tilde Q$ and $q$, one can fix
$\lambda^\#< \tilde \lambda$ to be the largest number with the
property $\tilde Q(\lambda^\#)=q(\lambda^\#)$. Thus, $q-\tilde Q$ is
non-negative on $[\lambda^\#, \tilde \lambda]$. We define for
$\lambda>0$ the function $z_\lambda:(0,\infty)\longrightarrow
\mathbb R$ by
$z_\lambda(\rho):=C^{-2}(n-2)^{-2}\rho^\frac{2}{p}+\lambda\rho$. By
relations \eqref{minko-2-egyenlotleseg} and
\eqref{minko-3-egyenloseg}, for every $\lambda>0$, we have
\[ z_\lambda(-\tilde Q'(\lambda))\leq \frac{n-1+a}{1-a}\tilde Q(\lambda), \qquad
 z_\lambda(- q'(\lambda))= \frac{n-1+a}{1-a} q(\lambda). \]
Since $z_\lambda$ is increasing, one has in particular that
\[ \tilde Q'(\lambda)-q'(\lambda)
 \geq z_\lambda^{-1}\left(\frac{n-1+a}{1-a} q(\lambda)\right)
 -z_\lambda^{-1}\left(\frac{n-1+a}{1-a} \tilde Q(\lambda)\right),
 \quad \lambda\in [\lambda^\#, \tilde \lambda]. \]
Taking into account that $z_\lambda^{-1}$ is increasing and $q \geq
\tilde Q$ on $[\lambda^\#, \tilde \lambda]$, the above inequality
implies that $(\tilde Q-q)'(\lambda)\geq 0$ for every $\lambda\in
[\lambda^\#, \tilde \lambda]$. In particular, we obtain $0>(\tilde
Q-q)(\tilde \lambda)\geq (\tilde Q-q)(\lambda^\#)=0$, a
contradiction. This completes the proof of
(\ref{g-kozott-egyenlotlenseg}).
\end{step}

\begin{step}\rm
We continue to assume $C>K_a$. Observe from \eqref{Q-definicio-3},
\eqref{Q-tilde-def} and \eqref{g-kozott-egyenlotlenseg} that
\begin{equation}\label{hasonlosag---o}
 \int_0^\infty \left\{ \mu(B(x_0,\rho))
 -(C^{-1}K_a)^\frac{n}{1-a}\mu_E(\mathbb B_n(\rho))\right\}
 f(\lambda,\rho) \,d\rho\geq 0,\quad \lambda>0.
\end{equation}
By the hypothesis $({\bf VD})^n_{C_0}$, for every $\rho>0$, we have
\[
C_0\frac{\mu(B(x_0,\rho))}{\mu_E(\mathbb B_n(\rho))}
 \geq \limsup_{r\to \infty} \frac{\mu(B(x_0,r))}{\mu_E(\mathbb B_n(r))}
 =:s_0.
\]
We claim that
\begin{equation}\label{s_0-becsles}
 s_0\geq (C^{-1}K_a)^\frac{n}{1-a}.
\end{equation}
Assuming the contrary, there exists $\delta_0>0$ such that, for some
$r_0>0$,
\[
\frac{\mu(B(x_0,\rho))}{\mu_E(\mathbb B_n(\rho))}
 \leq (C^{-1}K_a)^\frac{n}{1-a}-\delta_0
 \quad \text{for all}\ \rho\geq r_0.
\]
Hence, from \eqref{hasonlosag---o}, \eqref{two-measures} and
\eqref{Q-definicio-3}, we first have
\begin{align*}
0  &\le \int_0^\infty
 \left\{\mu(B(x_0,\rho))
 -(C^{-1}K_a)^\frac{n}{1-a} \mu_E(\mathbb B_n(\rho)) \right\}
 f(\lambda,\rho) \,d\rho \\
 &\le \int_0^{r_0} \mu(B(x_0,\rho)) f(\lambda,\rho) \,d\rho
 +\left\{ (C^{-1}K_a)^\frac{n}{1-a}-\delta_0\right\}
 \int_{r_0}^{\infty}\mu_E(\mathbb B_n(\rho))f(\lambda,\rho) \,d\rho \\
 &\quad -(C^{-1}K_a)^\frac{n}{1-a} \int_0^\infty
 \mu_E(\mathbb B_n(\rho))f(\lambda,\rho) \,d\rho\\
 &\le C_0\int_0^{r_0} \mu_E(\mathbb B_n(\rho))f(\lambda,\rho) \,d\rho
 -(C^{-1}K_a)^\frac{n}{1-a}
 \int_0^{r_0}\mu_E(\mathbb B_n(\rho))f(\lambda,\rho) \,d\rho\\
 &\quad -\delta_0\int_{r_0}^{\infty}\mu_E(\mathbb B_n(\rho))f(\lambda,\rho) \,d\rho\\
 &=\left\{ C_0-(C^{-1}K_a)^\frac{n}{1-a}+\delta_0\right\}
 \int_0^{r_0}\mu_E(\mathbb B_n(\rho))f(\lambda,\rho) \,d\rho
 -\delta_0\int_{0}^{\infty}\mu_E(\mathbb B_n(\rho))f(\lambda,\rho) \,d\rho\\
 &= \left\{ C_0-(C^{-1}K_a)^\frac{n}{1-a}+\delta_0\right\}
 \int_0^{r_0}\mu_E(\mathbb B_n(\rho))f(\lambda,\rho) \,d\rho
 -\delta_0\frac{n-1+a}{1-a}\lambda^\frac{n-2+2a}{2(a-1)}Q_E(1),
\end{align*}
where we used $Q_E(\lambda)=\lambda^{\frac{n-2+2a}{2(a-1)}}Q_E(1)$
following from \eqref{Q_E}. Next, by using the explicit form
\eqref{f-definic} of $f(\lambda,\rho)$ and $a-1<0$, the following
estimate holds:
\begin{align*}
\int_0^{r_0}\rho^n f(\lambda,\rho) \,d\rho &\le
\lambda^{\frac{n}{a-1}}\int_0^{r_0}\rho^{n-1-ap}
 \left\{ \rho^{2-ap}\left(\frac{n-1+a}{1-a}(2-ap)+ap\right)+ap\lambda\right\}
 d\rho \\
&= \left( \frac{n-1+a}{1-a}(2-ap)+ap\right)
 \frac{r_0^{n-2ap+2}}{n-2ap+2}\lambda^{\frac{n}{a-1}}
 +ap\frac{r_0^{n-ap}}{n-ap}\lambda^{\frac{n}{a-1}+1}.
\end{align*}
Reorganizing the above two estimates, we obtain the inequality of
type
\begin{equation}\label{M123}
M_1\lambda^\frac{n-2+2a}{2(a-1)}\leq
M_2\lambda^{\frac{n}{a-1}}+M_3\lambda^{\frac{n}{a-1}+1} \quad
\text{for all}\ \lambda>0,
\end{equation}
where $M_1,M_2,M_3>0$ are constants independent of $\lambda>0$.
Since
\[ \frac{n}{a-1}+1-\frac{n-2+2a}{2(a-1)}=\frac{n}{2(a-1)}<0, \]
\eqref{M123} fails for large values of $\lambda>0$. This
contradiction shows the validity of \eqref{s_0-becsles}.
\end{step}

\begin{step}\rm
Fix any $x\in X$. Since $B(x_0,r-d(x_0,x))\subset B(x,r)\subset
B(x_0,r+d(x_0,x))$ for every $r>d(x_0,x)$, on account of the
hypothesis $({\bf VD})^n_{C_0}$ and (\ref{s_0-becsles}), one has
\[ C_0\frac{\mu(B(x,\rho))}{\mu_E(\mathbb B_n(\rho))}
 \geq\limsup_{r\to \infty} \frac{\mu(B(x,r))}{\mu_E(\mathbb B_n(r))}
 =\limsup_{r\to \infty} \frac{\mu(B(x_0,r))}{\mu_E(\mathbb B_n(r))}
 =s_0\geq (C^{-1}K_a)^\frac{n}{1-a} \]
for all $\rho>0$. This concludes the proof. \hfill $\square$
\end{step}

\section{Applications: Caffarelli-Kohn-Nirenberg inequality on Finsler manifolds}\label{sect-3}

Before proving Theorems~\ref{theorem-Finsler} and
\ref{theorem-Berwald}, we concisely recall some notions from the
theory of Finsler manifolds (see Bao, Chern and Shen~\cite{BCS},
Shen~\cite{Shen-lec} and Ohta~\cite{Ohta-Fint} for details), and
prove  the validity of the Caffarelli-Kohn-Nirenberg inequality on
Minkowski spaces.

\subsection{Preliminary notions from Finsler geometry}

\subsubsection{Finsler manifolds}

Let $M$ be a connected $n$-dimensional $C^{\infty}$-manifold and
$TM=\bigcup_{x \in M}T_{x} M $ be its tangent bundle.

\begin{definition}[Finsler manifolds]\rm
The pair $(M,F)$ is called a \textit{Finsler manifold} if a
continuous function $F:TM\longrightarrow [0,\infty)$ satisfies the
conditions:
\begin{enumerate}[(1)]
\item $F\in C^{\infty}(TM\setminus \{0\})$;
\item $F(x,tv)=|t|F(x,v)$ for all $t\in \mathbb R$ and $(x,v)\in TM$;
\item the $n \times n$ matrix
\begin{equation}\label{g_v}
g_{ij}(x,v):=\frac{1}{2}\frac{\partial^2 (F^2)}{\partial v^{i}
\partial v^{j}}(x,v),
 \quad \text{where}\ v=\sum_{i=1}^n v^i \frac{\partial}{\partial x^i},
\end{equation}
is positive definite for all $(x,v)\in TM\setminus\{ 0 \}$. We will
denote by $g_v$ the inner product on $T_xM$ induced from
\eqref{g_v}.
\end{enumerate}
\end{definition}

If (and only if) $g_{ij}(x,v)$ is independent of $v$ in each $T_xM
\setminus \{0\}$, then $(M,F)$ gives a Riemannian manifold. A {\it
Minkowski space} consists of a finite dimensional vector space $V$
and a Minkowski norm which induces a Finsler metric on $V$ by
translation (i.e., $F(x,v)$ is independent of $x$). A Finsler
manifold $(M,F)$ is called a {\it locally Minkowski space} if any
point in $M$ admits a local coordinate system $(x^i)$ on its
neighborhood such that $F(x,v)$ depends only on $v$ and not on $x$.

For a $C^{\infty}$-curve $\sigma: [0,l]\longrightarrow M$, its
\textit{integral length} is given by $L_F(\sigma):=\int_{0}^{l}
F(\sigma, \dot\sigma)\,dt$. Define the {\it distance function}
$d_{F}: M\times M \longrightarrow[0,\infty)$ by $d_{F}(x_1,x_2) :=
\inf_{\sigma} L_F(\sigma)$, where $\sigma$ runs over all
$C^{\infty}$-curves from $x_1$ to $x_2$. When $(M,F)=(\mathbb
R^n,F)$ is a Minkowski space, one has $d_F(x_1,x_2)=F(x_2-x_1)$. A
$C^{\infty}$-curve $\sigma:[0,l] \longrightarrow M$ is called a
\emph{geodesic} if it is locally $d_F$-minimizing and has a constant
speed (i.e., $F(\sigma,\dot{\sigma})$ is constant). We can write
down the geodesic (Euler-Lagrange) equation in terms of the
\emph{covariant derivative} along $\sigma$ (see \cite{BCS} for
details). We say that $(M,F)$ is \emph{complete} if any geodesic
$\sigma:[0,l] \longrightarrow M$ can be extended to a geodesic
$\sigma:\mathbb R \longrightarrow M$.

The {\it polar transform} (or the \emph{dual norm}) of $F$ is
defined for every $(x,\alpha)\in T^*M$ by
\[ F^*(x,\alpha)
 :=\sup_{v\in T_xM\setminus \{0\}}\frac{\alpha(v)}{F(x,v)}. \]
Note that, for every $x\in M$, the function $F^*(x,\cdot)$ is a
Minkowski norm on $T_x^*M$. In particular, if $(\mathbb R^n,F)$ is a
Minkowski space, then so is $(\mathbb R^n,F^*)$ as well. For
$u(x)=d_F(x_0,x)$ with some fixed $x_0 \in M$, one can easily see
that $F^*(x,Du(x))=1$ for a.e.\ $x\in M$.

\subsubsection{Jacobi fields, Ricci curvature and volume comparison}

Let $\sigma:(-\varepsilon,\varepsilon) \times [0,l] \longrightarrow
M$ be a $C^{\infty}$-geodesic variation (i.e., $t \longmapsto
\sigma(s,t)$ is geodesic for each $s$), and put
$\eta(t)=\sigma(0,t)$. Then the variational vector field
$J(t):=\frac{\partial \sigma}{\partial s}(0,t)$ satisfies the
\emph{Jacobi equation}
\begin{equation}\label{Jacobi-egyenlet}
D^{\dot{\eta}}_{\dot{\eta}} D^{\dot{\eta}}_{\dot{\eta}}J
 +R^{\dot{\eta}}(J,\dot{\eta})\dot{\eta} \equiv 0,
\end{equation}
where $D^{\dot{\eta}}$ is the covariant derivative with reference
vector $\dot{\eta}$, and $R^{\dot{\eta}}$ is the curvature tensor
(see \cite{BCS} for details). For two linearly independent vectors
$v,w \in T_xM$ and $\mathcal S=\mathrm{span}\{v,w\}$, the {\it flag
curvature} of the \emph{flag} $(\mathcal S;v)$ is defined by
\[
K(\mathcal S;v) :=\frac{g_v(R^v(w,v)v,w)}{F(v)^2 g_v(w,w) -
g_v(v,w)^{2}}.
\]
If $(M,F)$ is Riemannian, then the flag curvature reduces to the
sectional curvature which depends only on $\mathcal S$ (not on the
choice of $v \in \mathcal S$). Take $v\in T_xM$ with $F(x,v)=1$ and
let $\{e_i\}_{i=1}^n$ with $e_n=v$ be an orthonormal basis of
$(T_xM,g_v)$ for $g_v$ from \eqref{g_v}. Put $\mathcal S_i={\rm
span}\{e_i,v\}$ for $i=1,...,n-1$. Then the {\it Ricci curvature} of
$v$ is defined by ${\rm Ric}(v):=\sum_{i=1}^{n-1}K(\mathcal S_i;v)$.
For $c \ge 0$, we also set $\mathrm{Ric}(cv):=c^2\mathrm{Ric}(v)$.

Shen gave a useful interpretation of these Finsler curvatures from
the Riemannian viewpoint (see \cite[\S 6.2]{Shen-lec}). Fix $v \in
T_xM \setminus \{0\}$ and extend it to a $C^{\infty}$-vector field
$V$ around $x$ such that all integral curves of $V$ are geodesic.
Then the flag curvature $K(\mathcal S;v)$ coincides with the
sectional curvature of $\mathcal S$ with respect to the Riemannian
structure $g_V$, and $\mathrm{Ric}(v)$ coincides with the Ricci
curvature of $v$ with respect to $g_V$. This observation leads the
following definition of the weighted Ricci curvature associated with
an arbitrary measure on $M$. We refer to \cite{Ohta-Fint},
\cite{Ohta-Sturm}, \cite{OS-BW}, and \cite{Ohta-Fspl} for details
and applications.

\begin{definition}[Weighted Ricci curvature]\rm
Let $\mu$ be a positive $C^{\infty}$-measure on $M$. Given $v \in
T_xM \setminus \{0\}$, let $\sigma:(-\varepsilon,\varepsilon)
\longrightarrow M$ be the geodesic with $\dot{\sigma}(0)=v$ and
decompose $\mu$ along $\sigma$ as
$\mu=e^{-\psi}\mathrm{vol}_{\dot{\sigma}}$, where
$\mathrm{vol}_{\dot{\sigma}}$ denotes the volume form of the
Riemannian structure $g_{\dot{\sigma}}$. Then, for $N \in
[n,\infty]$, the \emph{$N$-Ricci curvature} $\mathrm{Ric}_N$ is
defined by
\[ \mathrm{Ric}_N(v):=\mathrm{Ric}(v) +(\psi \circ \sigma)''(0)
 -\frac{(\psi \circ \sigma)'(0)^2}{N-n}, \]
where the third term is understood as $0$ if $N=\infty$ or if $N=n$
with $(\psi \circ \sigma)'(0)=0$, and as $-\infty$ if $N=n$ with
$(\psi \circ \sigma)'(0) \neq 0$.
\end{definition}

In particular, $\mathrm{Ric}_n$ is bounded below only if $(\psi
\circ \sigma)' \equiv 0$ along any $\sigma$. In terms of
$\mathrm{Ric}_N$, one can show the following Bishop-Gromov-type
volume comparison theorem. (Indeed, we can reduce it to the
Riemannian setting by using the gradient vector field of the
distance function from the center $x$.) We state only the
non-negatively curved case.

\begin{theorem} {\rm (\cite[Theorem~7.3]{Ohta-Fint})}\label{volume-comp-fo}
Let $(M,F,\mu)$ be a complete $n$-dimensional Finsler manifold with
non-negative $N$-Ricci curvature. Then we have
\[ \frac{\mu(B(x,R))}{\mu(B(x,r))} \le \left(\frac{R}{r}\right)^N
 \qquad \text{for every}\ x\in M,\ 0<r<R. \]
Moreover, if equality holds with $N=n$ for all $x \in M$ and
$0<r<R$, then any Jacobi field $J$ along a geodesic $\sigma$ has the
form $J(t)=tP(t)$, where $P$ is a parallel vector field along
$\sigma$ $($i.e., $D^{\dot{\sigma}}_{\dot{\sigma}}P \equiv 0)$.
\end{theorem}

We will actually need only the most restrictive case of $N=n$.

\subsection{Caffarelli-Kohn-Nirenberg inequality on Minkowski spaces}

Let $(M,F)$ be a Finsler manifold and $u\in {\rm Lip}_0(M)$. Note
that the local Lipschitz constant of $u$ is given by
$|Du|(x)=F^*(x,Du(x))$ for a.e.\ $x\in M$. Therefore, due to density
reasons, the Caffarelli-Kohn-Nirenberg inequality $({\bf
CKN})_C^{x_0}$ in the Finsler context takes the more familiar form
\[ \left(\int_M \frac{|u(x)|^p}{d_F(x_0,x)^{ap}} \,d\mu(x)\right)^\frac{1}{p}
 \leq C \left(\int_M {F^*}(x,D u(x))^2 \,d\mu(x)\right)^\frac{1}{2}
 \quad {\rm for\ all}\ u\in C_0^\infty(M). \]
We first prove  that the \emph{sharp} Caffarelli-Kohn-Nirenberg
inequality (with $C=K_a$ from \eqref{ckn-eukl}) holds on an
arbitrary Minkowski space $(\mathbb R^n,F)$ endowed with the
Lebesgue measure $\mu_F$ normalized so that
$\mu_F(B(0,1))=\omega_n$.

Let us first recall two inequalities on $(\mathbb R^n,F,\mu_F)$.
Given a measurable set $\Omega\subset \mathbb R^n$, let us denote by
$\Omega^\star$ the {\it anisotropic symmetrization} of $\Omega$,
i.e., it is the open ball with center $0$ such that
$\mu_F(\Omega^\star)=\mu_F(\Omega)$. For a function $u:\mathbb R^n
\longrightarrow \mathbb R$, $u^\star(x) :=\sup\{c\in \mathbb R:x\in
\{u>c\}^\star\}$ is the {\it anisotropic $($decreasing$)$
symmetrization} of $u$, where $\{u>c\}=\{x\in \mathbb R^n:u(x)>c\}$.
Due to Alvino,  Ferone,  Lions and Trombetti~\cite{AIHP-Lions} and
Van Schaftingen~\cite{vanSch}, one has
\begin{itemize}
\item {\it anisotropic P\'olya-Szeg\H{o} inequality}:
\[ \int_{\mathbb R^n} F^*(D u^\star(x))^2 \,d\mu_F(x)
 \leq \int_{\mathbb R^n} F^*(Du(x))^2 \,d\mu_F(x)
 \quad {\rm for\ all}\ u\in C_0^\infty(\mathbb R^n)_+; \]
\item {\it anisotropic Hardy-Littlewood inequality}:
if $p>1$ and $a\in [0,1]$, then we have
\[ \int_{\mathbb R^n} \frac{u(x)^p}{F(x)^{ap}} \,d\mu_F(x)
 \leq \int_{\mathbb R^n} \frac{u^\star(x)^p}{F(x)^{ap}} \,d\mu_F(x)
 \quad {\rm for\ all}\ u\in C_0^\infty(\mathbb R^n)_+, \]
where $C_0^\infty(\mathbb R^n)_+:=\{u\in C_0^\infty(\mathbb R^n):
u\geq 0\}$.
\end{itemize}

\begin{proposition}\label{prop-1-Minko}
Let $(\mathbb R^n, F)$ be a Minkowski space with $n\geq 3$, $x_0\in
\mathbb R^n$, $a\in [0,1)$, and $p={2n}/({n-2+2a})$. Then the sharp
Caffarelli-Kohn-Nirenberg inequality $({\bf CKN})_{K_a}^{x_0}$ holds
on $(\mathbb R^n, F,\mu_F)$. Moreover, the constant $K_a$ is optimal
and a family of extremals is given by
\[ u_\lambda(x)={\left(\lambda+F(x-x_0)^{2-ap}\right)^{\frac{2-n}{2-ap}}},
 \quad \lambda>0. \]
\end{proposition}

\noindent{\it Proof.} Recall that $d_F(x_0,x)=F(x-x_0)$. Without
loss of generality, we may assume that $x_0=0$. Let us consider the
constant
\[ C_a=\inf_{u\in C_0^\infty(\mathbb R^n)\setminus \{0\}}
 \frac{\left(\int_{\mathbb R^n} {F^*}(D u(x))^2 \,d\mu_F(x)\right)^{1/2}}
 {\left(\int_{\mathbb R^n} |u(x)|^p F(x)^{-ap} \,d\mu_F(x)\right)^{1/p}}. \]
We shall claim that $C_a=K_a^{-1}$. Due to the reversibility of $F$,
it is enough to consider non-negative functions in the above
expression. By the anisotropic P\'olya-Szeg\H{o} and
Hardy-Littlewood inequalities we have
\[ C_a=\inf_{u\in C_0^\infty(\mathbb R^n)_+\setminus \{0\}}
 \frac{\left(\int_{\mathbb R^n} {F^*}(D u^\star(x))^2 \,d\mu_F(x)\right)^{1/2}}
 {\left(\int_{\mathbb R^n} u^\star(x)^p F(x)^{-ap} \,d\mu_F(x)\right)^{1/p}}. \]
We may assume that $u^\star\in C_0^1(\mathbb R^n)_+$ (otherwise, a
density argument applies). Then there exists a non-increasing
function $h:[0,\infty)\longrightarrow [0,\infty)$ of class $C^1$
such that $u^\star(x)=h(F(x))$, and we have
\[ F^*(D u^\star(x))=F^*(h'(F(x))D F(x))
 =-h'(F(x))F^*(D F(x))=-h'(F(x)). \]
Therefore a simple calculation yields
\begin{equation}\label{Lieb-Talenti}
 \frac{\left(\int_{\mathbb R^n} {F^*}^2(D u^\star(x)) \,d\mu_F(x)\right)^{1/2}}
 {\left(\int_{\mathbb R^n} u^\star(x)^p F(x)^{-ap} \,d\mu_F(x)\right)^{1/p}}
 =\alpha_n^{\frac{1}{2}-\frac{1}{p}}
 \frac{\left(\int_0^\infty h'(\rho)^2 \rho^{n-1} \,d\rho\right)^{1/2}}
 {\left(\int_0^\infty h(\rho)^p \rho^{n-1-ap} \,d\rho\right)^{1/p}},
\end{equation}
where $\alpha_n=n\omega_n$ denotes the area of the unit sphere in
$\mathbb R^n$. On the other hand, following the approaches of
Lieb~\cite{Lieb} and Talenti~\cite{Talenti} in the Euclidean case
where the standard Schwarz symmetrization is used, one can see that
the minimizing expression is precisely the RHS of
\eqref{Lieb-Talenti}. Therefore, we have $C_a=K_a^{-1}$ which proves
our claim. Moreover, a class of minimizers $h_\lambda$ for
\eqref{Lieb-Talenti} is
$h_\lambda(\rho)={\left(\lambda+\rho^{2-ap}\right)^\frac{2-n}{2-ap}}$,
$ \lambda>0$, which can be obtained by the standard Euler-Lagrange
method. \hfill $\square$

\begin{remark}\rm
After a slight modification, Proposition~\ref{prop-1-Minko} remains
valid also for only \emph{positively homogeneous} Minkowski norms
(i.e., $F(tv)=tF(v)$ only for $t>0$). In such a case, the
anisotropic symmetrization is considered with respect to the
backward ball $B_-(0,1)=\{x\in \mathbb R^n:F(-x)< 1\}$, and the
level sets of the extremals have backward Wulff-shapes, homothetic
to $B_-(0,1)$ (see Krist\'aly \cite{Kristaly} and Van
Shaftingen~\cite{vanSch}).
\end{remark}

\subsection{Proof of Theorems~\ref{theorem-Finsler} and \ref{theorem-Berwald}}

\noindent{\it Proof of Theorem~$\ref{theorem-Finsler}$.} Since
$(M,F)$ is complete, by the Hopf-Rinow theorem it yields that
$(M,d_F,\mu)$ is a proper metric measure space. On account of
Theorem~\ref{volume-comp-fo}, $({\bf VD})^n_{C_0}$ holds with
$C_0=1$, while $\mu$ is normalized so as to satisfy $({\bf
AR})^n_{x_0}$. On the one hand, these properties imply that
\[ \mu(B(x,\rho)) \leq \mu_E(\mathbb B_n(\rho))
\quad {\rm for\ all}\ \rho>0,\ x\in M. \] On the other hand, by
$({\bf CKN})_{K_a}^{x_0}$, Theorem~\ref{theorem-1-fo} gives the
reverse inequality, thus equality holds. By
Theorem~\ref{volume-comp-fo}, it results that every Jacobi field $J$
along any geodesic $\sigma$ has the form $J(t)=tP(t)$, where $P$ is
a parallel vector field along $\sigma$. Then it follows from the
Jacobi equation \eqref{Jacobi-egyenlet} that
$R^{\dot\sigma}(J,\dot{\sigma})\dot{\sigma} \equiv 0$, so that
$K(\mathcal S;\dot\sigma) \equiv 0$ with $\mathcal S={\rm
span}\{\dot\sigma,P\}$. Due to the arbitrariness of $\sigma$ and
$J$, it turns out that the flag curvature of $(M,F)$ is identically
zero.
\hfill $\square$\\

\noindent{\it Proof of Theorem~$\ref{theorem-Berwald}$.} On the one
hand, since on every Berwald space $\mathrm{Ric}_n=\mathrm{Ric}$
holds for the Hausdorff measure $\mu_F$ (see
Shen~\cite[Propositions~2.6 \& 2.7]{Shen-volume}), we can apply
Theorem~\ref{theorem-Finsler} to see that the flag curvature of
$(M,F)$ is identically zero. On the other hand, every Berwald space
with zero flag curvature is necessarily a locally Minkowski space
(see Bao, Chern and Shen \cite[Section~10.5]{BCS}). Thanks to the
volume identity $\mu_F(B(x,\rho)) =\mu_E(\mathbb B_n(\rho))$,
$(M,F)$ must be isometric to a Minkowski space.
\hfill $\square$\\

We conclude the paper by presenting an example of a non-Riemannian
Berwald space.

\begin{example}\rm
We endow the space $\mathbb R^{n-1}$ ($n \geq 3$) with a Riemannian
metric $g$ such that $(\mathbb R^{n-1},g)$ is complete with
non-negative Ricci curvature. For every $\varepsilon> 0$, consider
on $\mathbb R^{n}=\mathbb R^{n-1}\times \mathbb R$ the metric
$F_\varepsilon:T\mathbb R^{n}\longrightarrow [0,\infty)$ given by
\[ F_\varepsilon((x,t),(v,w))=\sqrt{g_x(v,v)+w^2 +
\varepsilon \sqrt{g_x(v,v)^2+w^4}} \] for $(x,t)\in \mathbb R^{n}$,
$(v,w)\in T_{x}\mathbb R^{n-1}\times  T_t\mathbb R$. We observe that
$(\mathbb R^{n},F_\varepsilon)$ is a non-compact, complete,
non-Riemannian Berwald space with non-negative Ricci curvature.
According to Theorem~\ref{theorem-Berwald} and
Proposition~\ref{prop-1-Minko}, the following four statements are
equivalent:
\begin{itemize}
\item $({\bf CKN})_{K_a}^{\tilde x_0}$ holds on
$(\mathbb R^{n},F_\varepsilon,\mu_{F_\varepsilon})$ for {\it some}
element $\tilde x_0=(x_0,t_0)\in \mathbb R^{n}$;
\item $({\bf CKN})_{K_a}^{\tilde x_0}$ holds on
$(\mathbb R^{n},F_\varepsilon,\mu_{F_\varepsilon})$ for {\it every}
element $\tilde x_0=(x_0,t_0)\in \mathbb R^{n}$;
\item $g_x$ is independent of $x\in \mathbb R^{n-1}$
(i.e., $(\mathbb R^{n-1},g)$ is flat);
\item $(\mathbb R^{n},F_\varepsilon)$ is a Minkowski space.
\end{itemize}
\end{example}

\noindent {\bf Acknowledgment.} The paper was completed when A.
Krist\'aly visited the Department of Mathematics of Kyoto University
in October 2012. He is supported by a grant of the Romanian National
Authority for Scientific Research, CNCS-UEFISCDI, project number
PN-II-ID-PCE-2011-3-0241 and partially by Domus Hungarica. S. Ohta
is supported by the Grant-in-Aid for Young Scientists (B) 23740048.


\begin{thebibliography}{99}

\bibitem{AIHP-Lions} A. Alvino, V. Ferone, P.-L. Lions, G. Trombetti,
Convex symmetrization and applications. Ann. Inst. H. Poincar\'e
Anal. Non Lin\'eaire {\bf 14} (1997), no. 2, 275--293.

\bibitem{Aubin} T. Aubin, Probl\`emes isop\'erim\'etriques de Sobolev.
J. Differ. Geom. {\bf 11} (1976), 573--598.

\bibitem{BCS} D.~Bao, S.-S.~Chern, Z.~Shen, Introduction to Riemann--Finsler Geometry,
Graduate Texts in Mathematics, {\bf 200}, Springer Verlag, 2000.


\bibitem{CKN-cikk} L. Caffarelli, R. Kohn, L. Nirenberg, First order interpolation inequalities with weight.
Compos. Math. {\bf 53} (1984), 259--275.

\bibitem{CC-London} K. S. Chou, C. W. Chu,
On the best constant for a weighted Sobolev-Hardy inequality. J.
London Math. Soc. (2) {\bf 48} (1993), no. 1, 137--151.


\bibitem{doCarmo-Xia} M. P. do Carmo, C. Xia,
Complete manifolds with non-negative Ricci curvature and the
Caffarelli-Kohn-Nirenberg inequalities. Compos. Math. {\bf 140}
(2004), 818--826.


\bibitem{Heinonen} J. Heinonen, Lectures on Analysis on Metric Spaces,
Springer, New York, 2001.

\bibitem{Kristaly} A. Krist\'aly, Anisotropic singular phenomena in the presence of asymmetric
Minkowski norms. Preprint (2012).

\bibitem{Lieb} E. H. Lieb, Sharp constants in the Hardy-Littlewood and related inequalities.
Ann. of Math. (2) {\bf 118} (1983), 349--374.

\bibitem{Ohta-Fint} S. Ohta, Finsler interpolation inequalities.
Calc.\ Var.\ Partial Differential Equations {\bf 36} (2009),
211--249.

\bibitem{Ohta-Fspl} S.~Ohta,
Splitting theorems for Finsler manifolds of nonnegative Ricci curvature.
J.\ Reine Angew. Math.\ (to appear).
Available at {\sf arXiv:1203.0079}

\bibitem{Ohta-Sturm} S. Ohta, K.-T. Sturm, Heat flow on Finsler manifolds.
Comm. Pure Appl. Math. {\bf 62} (2009), no. 10, 1386--1433.

\bibitem{OS-BW} S. Ohta, K.-T. Sturm,
Bochner-Weitzenb\"ock formula and Li-Yau estimates on Finsler
manifolds. Preprint (2011). Available at {\sf arXiv:1104.5276}

\bibitem{Shen-volume} Z. Shen, Volume comparison and its applications in Riemann-Finsler geometry.
Adv. Math. {\bf 128} (1997), no. 2, 306--328.

\bibitem{Shen-lec} Z.~Shen, Lectures on Finsler geometry,
World Scientific Publishing Co., Singapore, 2001.

\bibitem{Szabo} Z. I. Szab\'o, Positive definite Berwald spaces. Structure theorems on Berwald spaces.
Tensor (N.S.) {\bf 35} (1981), no. 1, 25--39.

\bibitem{Talenti} G. Talenti, Best constant in Sobolev inequality.
 Ann. Mat. Pura Appl. (4) {\bf 110} (1976), 353--372.

\bibitem{vanSch} J. Van Schaftingen, Anisotropic symmetrization.
Ann. Inst. H. Poincar\'e Anal. Non Lin\'eaire {\bf 23} (2006),  no.
4, 539--565.

\end{thebibliography}
\end{document}